\documentclass[letterpaper, preprint, paper,11pt]{AAS}	

\usepackage{bm}
\usepackage{amsmath}
\usepackage{amssymb}
\usepackage{subcaption}
\usepackage[colorlinks=true, pdfstartview=FitV, linkcolor=black, citecolor= black, urlcolor= black]{hyperref}
\usepackage{overcite}
\usepackage{float}
\usepackage{footnpag}			      	

\PaperNumber{25-664}

\begin{document}

\title{Optimal Impulsive Control of Cislunar Relative Motion using Reachable Set Theory}

\author{Matthew Hunter\thanks{PhD Candidate, Department of Aeronautics and Astronautics, Stanford University, Stanford, CA 94305.**},  
\ Walter J. Manuel\thanks{PhD Candidate, Department of Aeronautics and Astronautics, Stanford University, Stanford, CA 94305.**},
\ and Simone D'Amico\thanks{Associate Professor, Department of Aeronautics and Astronautics, Stanford University, Stanford, CA 94305.\\ **The first two authors have contributed equally to the paper.}
}

\maketitle{}

\begin{abstract}
This work presents the first application of the state-of-the-art Koenig-D’Amico reachable set theory solver to cislunar, chaotic relative motion in the Circular-Restricted Three-Body Problem (CR3BP).
The relative motion dynamics of two spacecraft, a chief and a deputy, in the CR3BP are formulated as a Linear Time-Variant (LTV) system, allowing the solver to find an optimal impulsive control maneuver plan.
This methodology demonstrates robust and accurate control performance for both small and large reconfigurations over different CR3BP orbits and control windows.
These capabilities are enhanced by a Model Predictive Control (MPC) architecture to reject all sources of control, navigation, and dynamic error. 
The performance of the proposed approach is validated by unit testing, Monte Carlo simulations, and comparisons to baseline models for spacecraft relative motion.
Overall, this work demonstrates an optimal control methodology with the computational efficiency to be used on-board spacecraft, enabling the safe, effective, and efficient operation of Distributed Space Systems in cislunar space.
\end{abstract}

\section{Introduction}

The next era of space exploration is poised to be dually defined by activity in cislunar space and advancements in Distributed Space Systems (DSS).
The race to return to the moon is currently spearheaded by NASA's Artemis Program, which aims to land humans on the lunar surface by 2027 \cite{smith2020artemis}. 
NASA is also partnering with international allies, including the European Space Agency (ESA), Japan Aerospace Exploration Agency (JAXA), and Canadian Space Agency (CSA), to develop the Lunar Gateway station, which will be placed in a Near-Rectilinear Halo Orbit (NRHO) \cite{fuller2022gateway}.
Additionally, the defense sector, specifically the US Space Force, has identified increasing Space Domain Awareness (SDA) in cislunar space as a priority \cite{holzinger2021primer}.
This interest in cislunar space by government-focused sectors has already begun to drive investment from the commercial sector as well to help support and supplement civil and defense objectives.

Meanwhile, in Earth orbit, DSS are becoming increasingly commonplace, particularly in applications involving Rendezvous, Proximity Operations, and Docking (RPOD). 
DSS are defined as two or more spacecraft that work cooperatively, whether in formations, constellations, swarms, or other networks, to accomplish mission objectives not feasible with a single monolithic spacecraft. 
Multi-spacecraft systems will be an essential linchpin of the nascent cislunar ecosystem, as extensive RPOD operations of both crewed and uncrewed spacecraft are necessary to support lunar exploration and in-space logistics.

However, the cislunar regime features more complicated and chaotic dynamics compared to Earth orbit, motivating further work to ensure the safe, effective, and efficient operation of cislunar DSS. 
This requires innovation both in the dynamics modeling and optimal control strategies.
Dynamics modeling of two-body relative motion is fairly well explored, with widely-known approaches such as the Hill-Clohessy-Wiltshire equations\cite{clohessy_terminal_1960}, the Yamanaka-Ankersen state transition matrices\cite{yamanaka_new_2002}, and the D'Amico Relative Orbital Elements (ROE) \cite{damicophd}.
However, many of the assumptions used to create those models do not always hold in the cislunar environment, necessitating different dynamics models for three-body relative motion.
Recently, relative dynamics of the Circular Restricted Three-Body Problem (CR3BP) and the Elliptic Restricted Three-Body Problem (ER3BP) have been derived in the Local-Vertical-Local-Horizontal frame \cite{franzini_relative_2019, khoury_relative_2022}, as well as in the velocity-oriented TNW frame \cite{takubopassively} by linearizing around the chief spacecraft's position. 
Some authors have applied Linear Quadratic Regulator (LQR) control to a relative dynamics model obtained by linearizing the dynamics relative to a libration point, although this assumption potentially reduces accuracy in certain mission scenarios where the linearization is too coarse \cite{ardaens_control_2008}.
One optimal impulsive control method from Guffanti and D'Amico that leverages reachable set theory and primer vector theory \cite{guffanti2018integration}, has been applied to three-body absolute motion (linearized relative to a libration point)\cite{manuel2022optimal}, and three-body relative motion \cite{vela2025modal}.
However, these strategies all mostly focus on operational constraints, rather than provable optimality, and few offer the computational efficiency required to be implemented on-board a spacecraft.

Therefore, this work proposes a novel approach for modeling and controlling CR3BP relative spacecraft dynamics on-board over typical DSS control windows, resulting in three major contributions to the state of the art. 

\begin{enumerate}
    \item Two different methodologies are presented that leverage a linear instantaneous CR3BP relative dynamics model to generate a state transition matrix. The state transition matrix can be produced by either using a numerically efficient linear time-invariant approximation, breaking up the linear time-variant time horizon into smaller linear time invariant intervals, or using numerical integration, trading computational effort for dynamic accuracy. This facilitates the creation of a linear time-variant model of the relative CR3BP dynamics that propagates initial state and control actions to the end of a given time horizon.
    \item This modeling approach enables the first application of the state of the art Koenig-D'Amico solver, or ``KD solver," to the chaotic CR3BP dynamics of relative motion in the cislunar regime. The KD solver is a computationally-efficient, provably-optimal impulsive maneuver planner built upon a set of optimality conditions found through reachable set theory \cite{koenig_RST}. The KD solver has previously been applied to optimal path planning for satellites subject to two-body relative motion\cite{hunter_fast_2025}, and the same reachable set optimality conditions have been extended to derive closed form solutions for DSS \cite{chernick2021closed}.
    \item The differing approaches for CR3BP state transition matrix creation are demonstrated in extensive simulation analysis. Variations in solver optimality, computational effort, and dynamic accuracy are quantitatively evaluated both in isolation and comparatively through Monte Carlo analysis against well known two-body state transition matrices from the literature. A simple Model Predictive Control setup is also implemented to demonstrate dynamic error rejection under typical operational uncertainties.
\end{enumerate}

\section{Background}

\subsection{State Representation and Dynamics}

This section provides technical background on modeling the CR3BP dynamics of absolute and relative spacecraft motion. This paper considers a system comprised of two spacecraft, a chief and a deputy, operating in the vicinity of the Moon and influenced by the gravitational forces of both the Earth and the Moon. As a result, the CR3BP is an appropriate dynamics model to characterize the motion of the two spacecraft.

In the CR3BP, there are two bodies of significant mass, such as the Earth and the Moon, accompanied by a spacecraft of relatively negligible mass \cite{howell1984three}. 
The two larger bodies, known as the primaries, are modeled as point masses and rotate in circular orbits about their common center of mass. 
It is often convenient in CR3BP dynamic analysis to nondimensionalize the units of distance and time in the CR3BP and express motion in a synodic rotating frame that has the same inertial angular velocity as the primaries.

\begin{figure}[htb]
    \centering
    \begin{subfigure}[b]{0.4\textwidth}
        \centering
        \includegraphics[width=\textwidth]{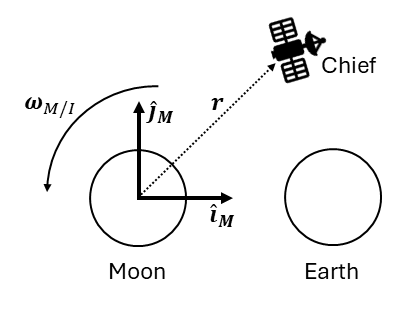}
        \caption{Moon reference frame}
        \label{fig:moon_frame}
    \end{subfigure}
    \begin{subfigure}[b]{0.25\textwidth}
        \centering
        \includegraphics[width=\textwidth]{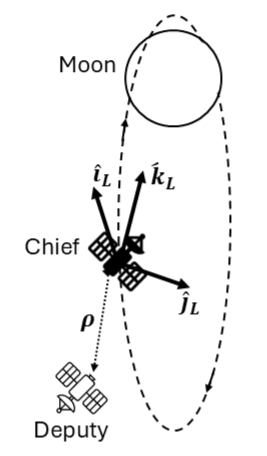}
        \caption{LVLH reference frame}
        \label{fig:lvlh_frame}
    \end{subfigure}  
    \caption{Key reference frames, vectors, and quantities used in this work. The chief spacecraft is depicted as being on a Near-Rectilinear Halo Orbit (NRHO).}
    \label{fig:frames} 
\end{figure}

The system of deputy and chief spacecraft is visualized in Figure \ref{fig:frames}, along with the two key reference frames employed in this paper. The chief spacecraft is depicted as being on a Near-Rectilinear Halo Orbit (NRHO). NRHOs are a subclass of halo orbits, three-dimensional periodic orbits that naturally form about the equilibrium points of the CR3BP system. Close proximity to the Moon and favorable stability properties distinguish NRHOs from other halo orbits \cite{zimovan2017characteristics}.

First, the moon-centered synodic rotating frame $\mathcal{M}$, as shown in Figure \ref{fig:moon_frame}, is defined by unit vectors $\hat{\bm{i}}_{M}$, $\hat{\bm{j}}_{M}$, and $\hat{\bm{k}}_{M}$, where $\hat{\bm{i}}_{M}$ points from the moon to the Earth, $\hat{\bm{k}}_{M}$ is aligned with the angular momentum vector of the Earth-Moon system, and $\hat{\bm{j}}_{M}$ completes the right-handed triad. 
Second, a Local-Vertical Local Horizontal (LVLH) frame $\mathcal{L}$, centered on the chief spacecraft and anchored to the moon, is also employed as shown in Figure \ref{fig:lvlh_frame}. Frame $\mathcal{L}$ is defined by unit vectors $\hat{\bm{i}}_{L}$, $\hat{\bm{j}}_{L}$, and $\hat{\bm{k}}_{L}$, where $\hat{\bm{k}}_{L}$ points from the chief to the Moon, $\hat{\bm{j}}_{L}$ is anti-parallel to the angular momentum vector of the chief, and $\hat{\bm{i}}_{L}$ completes the right-handed triad.
The LVLH frame is sometimes alternatively defined as the Radial/Tangential/Normal (RTN) frame, using a permutation of the same unit vectors, $\{ -\hat{\bm{k}}_L, \hat{\bm{i}}_L, -\hat{\bm{j}}_L \}$.
The absolute position of the chief spacecraft relative to the Moon’s center of mass is denoted as $\bm{r}$, and the relative position of the deputy spacecraft with respect to the chief is expressed as $\bm{\rho}$.

By linearizing the gravitational accelerations exerted on the deputy about the chief's position in the LVLH frame, the relative dynamics in the CR3BP system can be formulated as a Linear Time-Variant (LTV) system of equations in state space form as 
\begin{align}
    \bm{\dot{x}} = \bm{A}(t)\bm{x} + \bm{B}\bm{u}, \label{eq: state space eom}
\end{align}
where the state vector $\bm{x}$ is the Cartesian position and velocity of the deputy spacecraft in the LVLH frame, $\begin{bmatrix} \bm{\rho} & \bm{\dot{\rho}} \end{bmatrix}^{\top}$. \cite{franzini_relative_2019}
The plant matrix $\bm{A}(t)$ is defined as
\begin{align}
    \bm{A}(t) = 
    \begin{bmatrix}
        \bm{0}_{3}                &        \bm{I}_3 \\
        \bm{A_{\dot{\rho}\rho}}(t)  &  -2\bm{\Omega}_{L/I}(t)
    \end{bmatrix}. \label{eq: cr3bp A}
\end{align}
where $\boldsymbol{0}_3$ and $\boldsymbol{I}_3$ are 3x3 zero and identity matrices.
The matrix $\bm{A_{\dot{\rho}\rho}}$ concisely combines the terms in the lower left quadrant, defined as
\begin{align} \label{eq:linear_cr3bp_rel_dyn}
    \bm{A_{\dot{\rho}\rho}}(t) = &-[\bm{\dot{\Omega}}_{L/I}]^{\mathcal{L}} -\bm{\Omega}^{2}_{L/I}(t) 
    - \frac{\mu}{r^{3}}\left(\bm{I} - 3\frac{\bm{r}\bm{r}^\top}{r^2}\right)\\
    &- \frac{1 - \mu}{\| \bm{r} + \bm{r}_{em} \|^3}\left(\bm{I} - 3\frac{(\bm{r} + \bm{r}_{em})(\bm{r} + \bm{r}_{em})^\top}{{\| \bm{r} + \bm{r}_{em} \|}^2}\right), \notag
\end{align}
where $\bm{r}_{em}$ is the position vector of the earth relative to the moon, and $\mu$ is the mass parameter of the CR3BP system \cite{howell1984three}. The matrix $\bm{\Omega}_{L/I} \in \mathbb{R}^{3\times3}$ is a skew-symmetric matrix associated with the angular velocity vector of frame $\mathcal{L}$ relative to the inertial frame $\mathcal{I}$, $\bm{\omega}_{L/I}$. Similarly, the skew-symmetric matrix $[\bm{\dot{\Omega}}_{L/I}]^{\mathcal{H}}\in \mathbb{R}^{3\times3}$ associated with the angular acceleration vector $[\bm{\dot{\omega}}_{L/I}]^{\mathcal{L}}$ is also introduced. 
These two matrices are defined as
\begin{align}\label{skew_matrix}
    \bm{\Omega}_{L/I} = \begin{bmatrix}
        \quad0 & -\omega^{z}_{L/I} & \quad\omega^{y}_{L/I}\\
        \quad\omega^{z}_{L/I} & \quad0 & -\omega^{x}_{L/I}\\
        -\omega^{y}_{L/I} & \quad\omega^{x}_{L/I} & \quad0
    \end{bmatrix},\quad
    [\bm{\dot{\Omega}}_{L/I}]^{\mathcal{H}} = \begin{bmatrix}
        \quad0 & -\dot{\omega}^{z}_{L/I} & \quad\dot{\omega}^{y}_{L/I}\\
        \quad\dot{\omega}^{z}_{L/I} & \quad0 & -\dot{\omega}^{x}_{L/I}\\
        -\dot{\omega}^{y}_{L/I} & \quad\dot{\omega}^{x}_{L/I} & \quad0
    \end{bmatrix}.
\end{align}
The control matrix $\boldsymbol{B}$ transforms control actions $\boldsymbol{u}$ into instantaneous changes in state $\boldsymbol{x}$. Control actions are specified here to isolate changes in relative velocity, given as
\begin{equation}
    \boldsymbol{B} = 
    \begin{bmatrix}
       \boldsymbol{0}_3 \\ \boldsymbol{I}_3 
    \end{bmatrix}.  \label{eq: control matrix}
\end{equation}
It is worth noting here that this work will use this LTV state space model both as the foundational model leveraged to create the proposed approach and as the reference ground truth of the validation simulation in the \textit{Results} section below.

\subsection{Problem Definition}

This paper addresses the impulsive or instantaneous control of linear time-varying dynamic systems, as defined by the following. 
Consider a generic LTV dynamics system, given as
\begin{equation}
    \boldsymbol{x}_f = \boldsymbol{\Phi}(t_0)\boldsymbol{x}_0 + \int_{t_0}^{t_f} \boldsymbol{\Phi}(\tau)\boldsymbol{B}(\tau)\boldsymbol{u}(\tau)d\tau,
\end{equation}
where $\boldsymbol{x}_0 \in \mathbb{R}^n$ is the initial state at time $t_0$, $\boldsymbol{x}_f \in \mathbb{R}^n$ is the final state at time $t_f$, $\boldsymbol{\Phi}(t) = \boldsymbol{\Phi}(t,t_f) \in \mathbb{R}^{n \times n}$ is the State Transition Matrix (STM) that propagates a given state from time $t$ to $t_f$, and $\boldsymbol{B}(t) \in \mathbb{R}^{n \times m}$ is the control matrix that converts control actions $\boldsymbol{u}(t) \in \mathbb{R}^m$ to changes in state $\boldsymbol{x}$. 
For simplicity, the pseudostate $\boldsymbol{\omega} \in \mathbb{R}^n$, the dynamic invariants of the control problem, and control matrix $\boldsymbol{\Gamma}(t) \in \mathbb{R}^{n \times m}$, the translation from control actions to changes in final state, are defined as,
\begin{equation}
    \boldsymbol{\omega} = \boldsymbol{x}_f - \boldsymbol{\Phi}(t_0)\boldsymbol{x}_0, \hspace{5mm} \boldsymbol{\Gamma}(t) = \boldsymbol{\Phi}(t)\boldsymbol{B}(t).
\end{equation}
Therefore, the optimal control problem (OCP) considered by this paper is defined as
\begin{align}
    \text{minimize: }& \int_{t_0}^{t_f} f(\boldsymbol{u}(\tau),\tau)d\tau \hspace{5mm} \text{subject to: } \boldsymbol{\omega} = \int_{t_0}^{t_f} \boldsymbol{\Gamma}(\tau)\boldsymbol{u}(\tau)d\tau. \label{eq: continuous problem}
\end{align}
Cost function $f$ is norm-like, such that it satisfies three properties: 1) $f$ is defined for all maneuvers occurring within $t_{0} \leq t \leq t_{f}$, 2) all sublevel sets of $f(\boldsymbol{u}(t),t)$ are convex and compact, and 3) $f(\alpha\boldsymbol{u}(t),t) = \alpha f(\boldsymbol{u}(t),t)$ for any $\alpha \geq 0$. 
Examples of norm-like cost functions and their corresponding spacecraft propulsion constraints include the L2-norm $||\boldsymbol{u}||_2$, where the spacecraft aligns a single thruster with the desired maneuver direction, and the L1-norm $||\boldsymbol{u}||_1$, where the spacecraft has a fixed attitude and three pairs of thrusters mounted on opposite sides. 
It is worth noting that the choice of a single cost function, rather than the more-general piecewise definition in previous works, is made for simplicity and does not effect the methodology detailed in later sections.

This work focuses on impulsive control actions, or instantaneous changes in state, assuming that the time required to complete each maneuver is short with respect to the secular and long-term periodic effects present in $\boldsymbol{\Phi}(t)$. 
Given there are no restrictive magnitude constraints on $\boldsymbol{u}$, time-varying control actions can be expressed as a set of $k$ impulses, and the unconstrained control problem considered by the remainder of this paper is given as
\begin{equation}
	\text{minimize: } c = \sum_{j=1}^{k} f_j(\boldsymbol{u}_j,t_j) \hspace{3mm} \text{subject to: } \boldsymbol{\omega} = \sum_{j=1}^{k} \boldsymbol{\Gamma}(t_j)\boldsymbol{u}_j, \label{eq: impulsive problem}
\end{equation}
where all $t_0 \leq t_j \leq t_f$ and $c$ is the minimum cost objective.

\subsection{Reachable Set Theory}

The control solver proposed by this work is built upon the optimality conditions found through reachable set theory, and this analysis technique is detailed below. Let $U(c,t)$ be the set of control actions $\boldsymbol{u}$ with a cost no greater than $c$ at time $t$, defined as
\begin{equation}
    U(c,t) = \{\boldsymbol{u}: f(\boldsymbol{u},t) \leq c \}. \label{eq: U}
\end{equation}
Let $S(c,t)$ be the set of all $\boldsymbol{\omega}$ that can be reached by a single maneuver $\boldsymbol{u}$ of cost no greater than $c$ at time $t$, given as
\begin{equation}
    S(c,t) = \{\boldsymbol{y}: \boldsymbol{y} = \boldsymbol{\Gamma}(t)\boldsymbol{u}, \boldsymbol{u} \in U(c,t) \}. \label{eq: S(c,t)}
\end{equation}
This set is not generally convex. Let $T$ be the discontinuous set of times $t$ for admissible control actions in a control window bounded by $t_0$ and $t_f$, given as
\begin{equation}
    T = \{t: t_0 \leq t \leq t_f \}. \label{eq: T}
\end{equation}
Then, the set of all $\boldsymbol{\omega}$ that can be reached by a single maneuver $\boldsymbol{u}$ of cost no greater than $c$ at any any time in $T$ is defined as
\begin{equation}
    S(c,T) = \cup_{t \in T} S(c,t), \label{eq: S(c,T)}
\end{equation}
which is also generally non-convex. Finally, let $S^*(c,T)$ be the set of $\boldsymbol{\omega}$ that can be reached by $k \geq 1$ control actions of combined cost no greater than $c$ at times in $T$, defined equivalently as
\begin{align}
    S^*(c,T) = \{\boldsymbol{z}: \boldsymbol{z} &= \sum\nolimits_{j=1}^{k} \boldsymbol{\Gamma}(t_j)\boldsymbol{u}_j, t_j \in T, \boldsymbol{u}_j \in U(c_j,t_j), \sum\nolimits_{j=1}^{k} c_j = c \} \label{eq: S* base} \\
    = \{\boldsymbol{z}: \boldsymbol{z} &= \sum\nolimits_{j=1}^{k}\alpha_j\boldsymbol{y}_j, \boldsymbol{y}_j \in S(c,T), \alpha_j \geq 0, \sum\nolimits_{j=1}^{k} \alpha_j = 1 \}. \label{eq: S* convex hull}
\end{align}
This shows that $S^*(c,T)$ is a linear combination and convex hull of $\boldsymbol{\omega}$ in $S(c,T)$. Returning to the optimal control problem in Eq. (\ref{eq: impulsive problem}), the set $S^*(c,T)$ in Eq. (\ref{eq: S* convex hull}) completely describes the dynamics constraint such that the problem can be defined as
\begin{equation}
    \text{minimize: } c \hspace{5mm} \text{subject to: } \boldsymbol{\omega} \in S^*(c,T). \label{eq: problem with RST}
\end{equation}

\subsection{Computationally-Efficient Control Solver}

\begin{figure} [h!]
    \centering
    \includegraphics[width=0.8\linewidth]{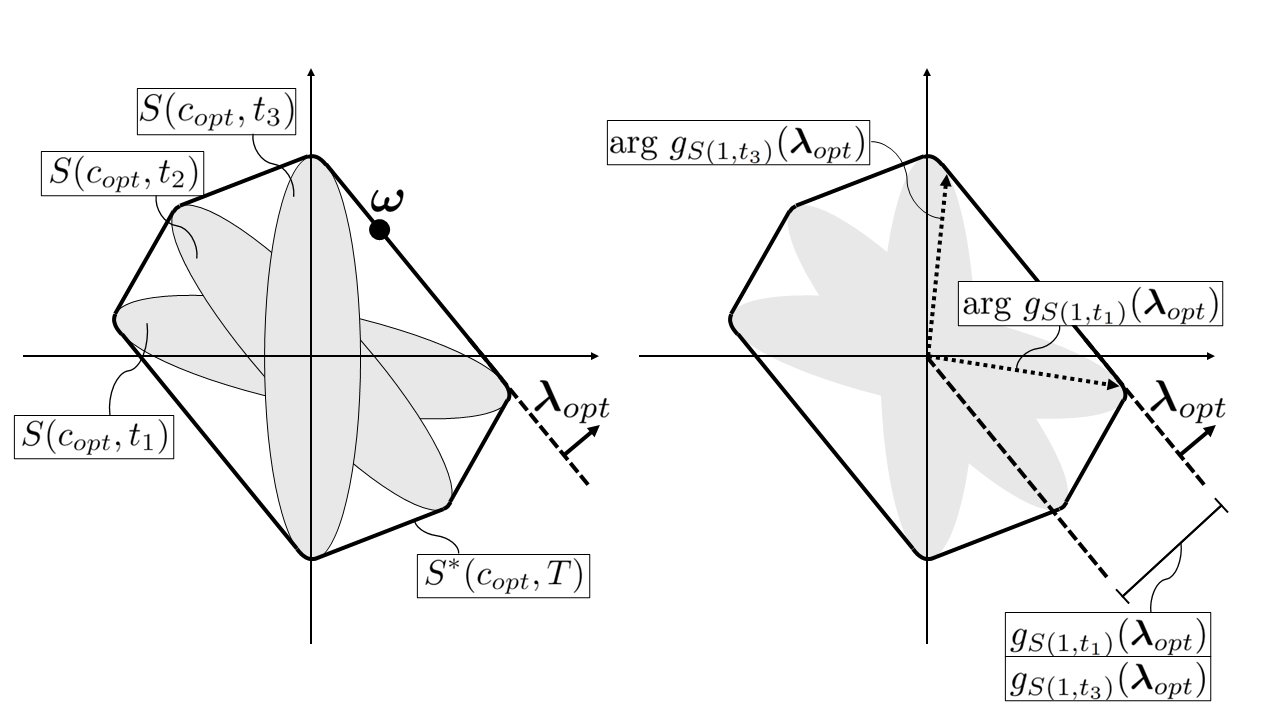}
    \caption{2-dimensional reachable set formed as convex hull of states reached by single maneuvers at times $t_1 < t_2 < t_3$ (left). Optimality conditions evaluated on the boundary of $S^*(1,T)$ through contact functions (right). \cite{hunter_fast_2025} }
    \label{fig: unconstrained}
\end{figure}

This work extends a computationally-efficient solver designed by Koenig and D'Amico, denoted in this work as the ``KD solver," for CR3BP application \cite{koenig_RST}. 
The KD solver is built upon a set of optimality conditions for control found through the Second-Order Cone Program (SOCP) dual of Eq. (\ref{eq: problem with RST}), given as 
\begin{equation}
    \text{maximize: } \boldsymbol{\lambda}^T\boldsymbol{\omega} \hspace{5mm} \text{subject to: } \max_{t \in T} g_{U(1,t)}(\boldsymbol{\Gamma}^T(t) \boldsymbol{\lambda}) \leq 1. \label{eq: RST SOCP}
\end{equation}
\noindent where contact function $g_{U(1,t)}(\boldsymbol{\Gamma}^T(t) \boldsymbol{\lambda}) = \max_{\boldsymbol{\boldsymbol{u}} \in U(1,t)} \boldsymbol{\lambda}^T\boldsymbol{\Gamma}(t)\boldsymbol{u}$.
First, consider an optimal solution to Eq. (\ref{eq: problem with RST}) and (\ref{eq: RST SOCP}) defined as $c_{opt}$ and $\boldsymbol{\lambda}_{opt}$. $c_{opt}$ is the minimum cost to reach $\boldsymbol{\omega}$ if and only if $\boldsymbol{\omega}$ is on the boundary of $S^*(c_{opt},T)$. 
This boundary is described by the supporting hyperplane of $S^*(c_{opt},T)$ defined perpendicular to $\boldsymbol{\lambda}_{opt}$, and the optimal control input profile may not be unique if $\boldsymbol{\omega}$ lies on an edge or point of the boundary of $S^*(c_{opt},T)$. 
Second, for any $\boldsymbol{\omega}$ on the boundary of $S^*(c_{opt},T)$, an optimal control input sequence must exist such that $k \leq n$.
Therefore, such a sequence must exist for any $\boldsymbol{\omega}$ of no more maneuvers than the dimensionality of $\boldsymbol{\omega}$. 
Third, control actions should only occur when $g_{U(1,t)}(\boldsymbol{\Gamma}^T(t) \boldsymbol{\lambda}_{opt})$ reaches its maximum value in the control window $T$. 
This condition is reflected in the constraint of Eq. (\ref{eq: RST SOCP}). 
Fourth, control actions should be in the direction of $\text{arg }g_{U(1,t)}(\boldsymbol{\Gamma}^T(t)\boldsymbol{\lambda}_{opt})$. 
The third and fourth conditions can be evaluated for an arbitrary $\boldsymbol{\lambda}$ and time $t$ for an L2-norm cost function through the Cauchy Schwartz inequality, given as
\begin{align}
    g_{U(1,t)}(\boldsymbol{\Gamma}^T(t) \boldsymbol{\lambda}) &= ||\boldsymbol{\Gamma}^T(t) \boldsymbol{\lambda}||_2, \label{eq: g_U} \\
    \text{arg }g_{U(1,t)}(\boldsymbol{\Gamma}^T(t) \boldsymbol{\lambda}) &= \frac{\boldsymbol{\Gamma}^T(t) \boldsymbol{\lambda}}{||\boldsymbol{\Gamma}^T(t) \boldsymbol{\lambda}||_2}. \label{eq: arg g_U}
\end{align}
These optimality conditions are visualized in Figure \ref{fig: unconstrained}.

The KD solver consists of three major components: Initialization, Refinement, and Control Input Extraction. 
Initialization uses the direction of $\boldsymbol{\omega}$ as a good initial guess of the solution $\boldsymbol{\lambda}_{opt}$ of Eq. (\ref{eq: RST SOCP}) to produce a set of initial candidate maneuvers times $T_{est}$, defined as the times where Eq. (\ref{eq: g_U}) reaches a large value relative to other times in $T$. 
Refinement iteratively solves Eq. (\ref{eq: RST SOCP}) to generate a new $\boldsymbol{\lambda}_{opt}$ for the times in $T_{est}$. 
This new $\boldsymbol{\lambda}_{opt}$ is used to remove sub-optimal candidate maneuver times from $T_{est}$ whose contact function Eq. (\ref{eq: g_U}) falls below a threshold $\epsilon_{remove}$, given as
\begin{equation}
    g_{U(1,t)}(\boldsymbol{\Gamma}^T(t)\boldsymbol{\lambda}_{opt}) < 1 - \epsilon_{remove}, \label{eq: eps remove threshold}
\end{equation}
and add optimal candidate maneuvers times from $T$ to $T_{est}$ that exceed an optimality tolerance $\epsilon_{cost}$, given as
\begin{equation}
    g_{U(1,t)}(\boldsymbol{\Gamma}^T(t)\boldsymbol{\lambda}_{opt}) > 1 + \epsilon_{cost}. \label{eq: eps cost threshold}
\end{equation}
Refinement is repeated until $g_{U(1,t)}(\boldsymbol{\Gamma}^T(t)\boldsymbol{\lambda}_{opt}) < 1 + \epsilon_{cost}$ for all $t \in T_{est}$. 
Control Input Extraction finds an unconstrained optimal maneuver plan, under only dynamic constraints, at the optimal maneuver times in $T_{est}$ by computing the optimal maneuver directions from Eq. (\ref{eq: arg g_U}) and finding the maneuver magnitudes that minimize the resulting pseudostate error of the maneuver plan through a Quadratic Program (QP). 
Limiting control inputs to times in $T_{est}$ rather than $T$ ensures that maneuvers only occur at optimal times without an explicit definition of optimality in the QP.

\section{Methodology}

To apply the KD solver to the CR3BP relative motion problem, a control matrix $\boldsymbol{\Gamma}$, comprised of an STM $\boldsymbol{\Phi}$ and control matrix $\boldsymbol{B}$, must be derived to form the dynamics constraint present in Eq. (\ref{eq: impulsive problem}). 
The control matrix $\boldsymbol{B}$ is simply chosen as the same used in the CR3BP state space equations in Eq. (\ref{eq: control matrix}).
Two different strategies are proposed by this work to create an STM from the CR3BP dynamics in Eq. (\ref{eq: cr3bp A}): matrix exponential and numerical integration.
For matrix exponential, a Linear Time-Invariant (LTI) STM $\boldsymbol{\Phi}_{LTI}$ can be produced from a dynamic plant matrix $\boldsymbol{A}$, given as
\begin{equation}
    \boldsymbol{\Phi}_{LTI}(t,t_f) = e^{\boldsymbol{A}(t)(t_f-t)}. \label{eq: LTI STM}
\end{equation}
Consider approximating the LTV dynamics as a series of LTI dynamic systems. 
By discretizing the control window into an arbitrary number of time steps $l$, the LTV STM $\boldsymbol{\Phi}$ can be defined as a multiplicative sequence of LTI STMs, given as
\begin{equation}
    \boldsymbol{\Phi}(t) = \Pi_{i=1}^{l-1} \boldsymbol{\Phi}_{LTI}(t_i,t_{i+1}),
\end{equation}
\noindent where $t_l = t_f$.
The accuracy of this approximation depends on both the instantaneous variance of the true nonlinear dynamic system and the density of the time discretization, resulting in a tradeoff between computational efficiency and dynamic accuracy.

For numerical integration, the LTV STM is found by integrating the linear instantaneous dynamics in plant matrix $\boldsymbol{A}$ from a given time $t$ to the end of the control window $t_f$. The time derivative of the STM $\boldsymbol{\Phi}$ is defined as
\begin{equation}
    \frac{\text{d}}{\text{d}t} \boldsymbol{\Phi}(t) = \boldsymbol{A}(t)\boldsymbol{\Phi}(t),\qquad
    \boldsymbol{\Phi}(t = t_0) = \boldsymbol{I}
\end{equation}
Therefore, the numerical integration problem resulting in the LTV STM is given as
\begin{equation}
    \boldsymbol{\Phi}(t) = \int_t^{t_f} \boldsymbol{A}(\tau)\boldsymbol{\Phi}(t,\tau-t) d\tau.
\end{equation}
The initial condition of this integration is the identity matrix, which logically follows from the definition of an STM over zero elapsed time. Numerical integration captures the full time-variance of the CR3BP relative dynamics over the control window.

\begin{figure}[h!]
    \centering
    \includegraphics[width=1\linewidth]{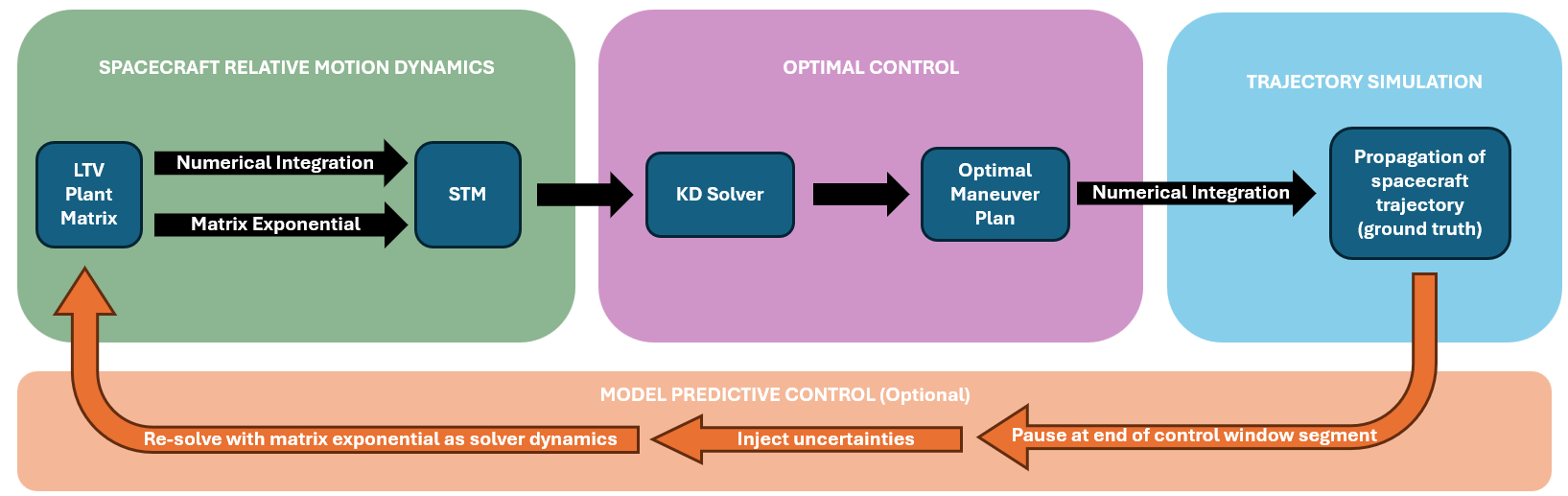}
    \caption{Block diagram depicting architecture of proposed methodology, from STM creation via two possible approaches, to solution of the OCP via the KD solver, to the implementation of that solution in numerically integrated ground truth dynamics. The option to close the loop with MPC control using the matrix exponential approach is also depicted. Other parameters not shown include the chief spacecraft absolute state, which is an input to form the LTV plant matrix; and inputs to the OCP inputs such as initial and final desired conditions, and control window length. }
    \label{fig: block}
\end{figure}

In summary, the two different strategies proposed here to create an STM from the CR3BP dynamics will demonstrate a tradeoff between computational performance and dynamic accuracy. Using closed-form matrix exponential solutions enables the LTV system to be efficiently approximated as a series of LTI systems. However, this approximation will lose accuracy in a highly variable dynamic environment and may need a prohibitively dense time discretization to recover usable control solutions. In contrast, numerical integration completely characterizes the CR3BP relative dynamics within the accuracy of the linear model, but requires more substantial computational resources to achieve this performance. Figure \ref{fig: block} outlines the overall open and closed loop architectures employed in this work, illustrating how an STM is generated from the LTV system through either of the proposed methods, and then that STM is used by the KD solver to solve for the optimal maneuver plan.

\section{Results}

The two proposed CR3BP STM creation approaches are thoroughly tested with the KD solver through multiple methods of validation. 
First, \textit{Initial Validation} consists of two relevant test cases for CR3BP DSS control to demonstrate and compare the capabilities of the matrix exponential and numerical integration approaches. 
Then, \textit{Validation Against Two-Body Dynamics Models} incorporates several commonly-used two-body STMs from the literature that have demonstrated a high level of accuracy for Earth orbit to show that the matrix exponential and numerical integration approaches outperform a two-body approximation of CR3BP motion. 
Next, \textit{Monte Carlo Validation} expands this comparative analysis to a series of randomized reconfigurations, to show that the conclusions from the test case analysis remain consistent over a variety of solver inputs. 
Finally, \textit{Model Predictive Control} incorporates the matrix exponential approach and KD solver into a simple robust control architecture, resolving the maneuver plan over an extended control window to reject dynamic modeling errors. 
For all validation simulations, control solutions are produced ``on-board" using the STM creation approach and corresponding LTV dynamics model. 
These maneuver plans are then numerically propagated by the unperturbed equations of CR3BP motion, used as the ground truth.
It is worth noting that this ground truth uses the same equations as that used to produce the matrix exponential and numerical integration STMs.
However, inherent dynamic differences from the ground truth are produced in both cases.
Ground truth errors in the matrix exponential approach are clearly created from the computationally efficient LTI approximation, while differing integration timesteps between the numerical integration approach and the simulation environment create small errors that are exaggerated by the chaotic dynamic environment around perilune.
Maneuver plan cost, final state error, and solver runtime will be primarily used to quantify control performance under different STM creation approaches. 
All simulations are conducted on the same hardware, using a 3.80 GHz processor and 32GB RAM.

\subsection{Initial Validation}

The proposed solver and approaches to formulate the LTV control problem are first validated against two general reconfigurations in cislunar orbit. The initial absolute orbits of the chief spacecraft for both reconfigurations are given in the moon-centered synodic frame in Table \ref{tab: initial chief states}.

\begin{table}[h!]
    \centering
    \begin{tabular}{c|c|c|c|c|c|c}
    Scenario  & $x$ (km) & $y$ (km) & $z$ (km) & $\dot{x}$ (km/s) & $\dot{y}$ (km/s) & $\dot{z}$ (km/s) \\
    \hline \hline
    Reconfiguration 1 & -13395 & 0 & -70841 & 0 & 0.1055 & 0 \\
    Reconfiguration 2 &  -4909 & 29088 & -14638 & 0.1080 & -0.1647 & 0.4331\\
    \hline
    \end{tabular}
    \caption{Absolute states of the chief spacecraft at $t_0$ in the moon-centered synodic frame for initial validation reconfigurations. }
    \label{tab: initial chief states}
\end{table}
\noindent Reconfiguration 1 represents a chief spacecraft in the 9:2 resonant NRHO at apolune, while Reconfiguration 2 represents a chief spacecraft in a 3:1 resonant halo orbit about to cross over perilune. 
The initial and final deputy relative states in LVLH frame for both reconfigurations are given in Table \ref{tab: initial and final deputy states}.
\begin{table}[]
    \centering
    \begin{tabular}{c|c|c|c|c|c|c}
    Deputy State  & $x$ (km) & $y$ (km) & $z$ (km) & $\dot{x}$ (km/s) & $\dot{y}$ (km/s) & $\dot{z}$ (km/s) \\
    \hline \hline
    \multicolumn{7}{c}{Reconfiguration 1} \\
    \hline
    Initial & -300 & -400 & -200 & 0 & 0 & 0 \\
    Final &  300 & 400 & 200 & 0 & 0 & 0 \\
    \hline
    \multicolumn{7}{c}{Reconfiguration 2} \\
    \hline
    Initial & -10 & -0.3 & -0.05 & 0 & 0 & 0 \\
    Final &  0.1 & 0.3 & 0.05 & 0 & 0 & 0 \\
    \hline
    \end{tabular}
    \caption{Initial and final deputy relative states in LVLH frame for initial validation reconfigurations.}
    \label{tab: initial and final deputy states}
\end{table}
\noindent The two initial reconfiguration scenarios were designed intentionally to represent plausible hypothetical applications in cislunar space. 
Reconfiguration 1 changes the leader-follower order of the chief and deputy spacecraft, such that the deputy begins by ``following" the chief with negative separation and ends by ``leading" the chief with positive separation.
Reconfiguration 1 has a 66.84h control window.
Reconfiguration 2 involves a reduction in inter-spacecraft separation to conduct proximity operations and docking near the chief spacecraft and has a 33.52h control window.
These reconfigurations will demonstrate the proposed solver's ability to accurately control CR3BP relative motion in real-world scenarios. 
An optimal control sequence for each reconfiguration is solved using both the matrix exponential and numerical integration STM approaches to compare resulting $\Delta v$ costs, final state errors, and computational efficiency.
These maneuver plans are simulated in a ground truth numerically integrated propagation of the unperturbed chief (nonlinear absolute CR3BP) and deputy (linear relative CR3BP) dynamics.

\begin{figure}[h!]
    \centering
    \includegraphics[width=1\linewidth]{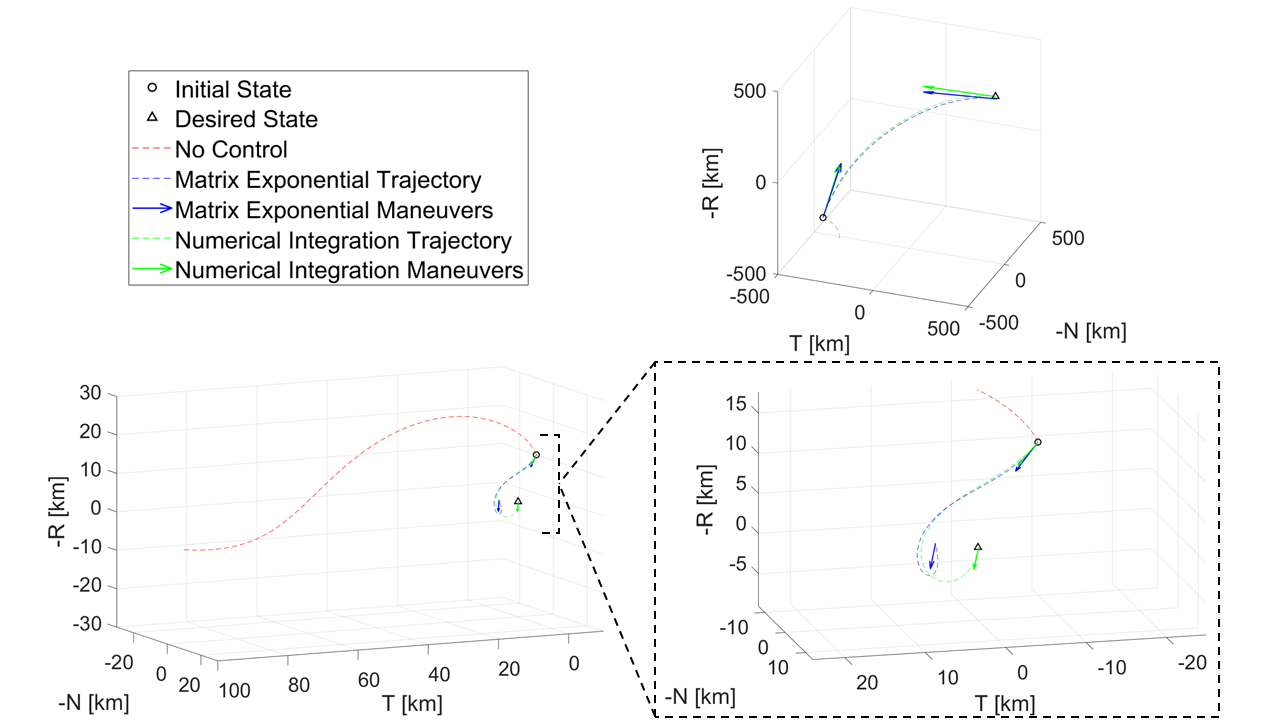}
    \caption{LVLH trajectories for Reconfiguration 1 (top row) and Reconfiguration 2 (bottom row).}
    \label{fig: reconfigs}
\end{figure}

The KD solver is set up with a control domain $T$ of 1001 candidate maneuver times for each reconfiguration, produced by uniformly splitting the control window into 1000 steps. 
The STM is precomputed at each candidate maneuver time before running the solver to isolate the computational effort associated with calculating the STMs from that required to find the optimal maneuver plan.
The matrix exponential approach breaks the Reconfiguration 1 control window into 10min time steps and the Reconfiguration 2 control window into 20min time steps for the LTI approximation.
The chosen cost function is the sum of L2-norms of RTN control actions at times in $T$, and control is admissible at all times in the control window.
The initialization component of the KD solver is set to evaluate the contact function for every ten candidate times in $T$, taking the ten highest as the initial $T_{est}$, and the Control Input Extraction component uses an identity matrix to minimize the pseudostate error.
The resulting control trajectories are displayed in Figure \ref{fig: reconfigs}.

For both reconfigurations, the KD solver found almost identical maneuver plans for the matrix exponential and numerical integration STM approaches, illustrated by the overlay of trajectories and maneuvers. 
The chief's proximity to perilune plays a clear role in the resulting controlled and uncontrolled trajectories.
The chief and deputy far from perilune in Reconfiguration 1 feature smooth, consistent changes in relative state with low curvature. 
Reconfiguration 2 demonstrates the highly variable dynamics near perilune with significant, high curvature changes in the deputy relative state path over a shorter control window than Reconfiguration 1.
The contrasting advantages and disadvantages of the two STM approaches are well illustrated by the performance metrics associated with each reconfiguration in Table \ref{tab: reconfigs}.
\begin{table}[h!]
    \centering
    \begin{tabular}{c|c|c|c|c}
        STM Strategy & Cost (m/s) & RMS $\boldsymbol{x}_f$ Error (km) & STM Runtime (s) & Solver Runtime (s) \\
        \hline \hline
        \multicolumn{5}{c}{Reconfiguration 1} \\
        \hline
        Matrix Exp. & 9.7644e-4 & 8.4613 & 0.55 & 0.96 \\
        Numerical Int. & 9.9372e-4 & 0.8065 & 8.63 & 1.00 \\
        \hline
        \multicolumn{5}{c}{Reconfiguration 2} \\
        \hline
        Matrix Exp. & 4.6352e-4 & 2.9950 & 13.36 & 26.84 \\
        Numerical Int. & 4.8193e-4 & 0.0496 & 26.17 & 26.40 \\
        \hline
    \end{tabular}
    \caption{Performance metrics for initial validation reconfigurations.}
    \label{tab: reconfigs}
\end{table}

\begin{figure}[h]
    \centering
    \includegraphics[width=1\linewidth]{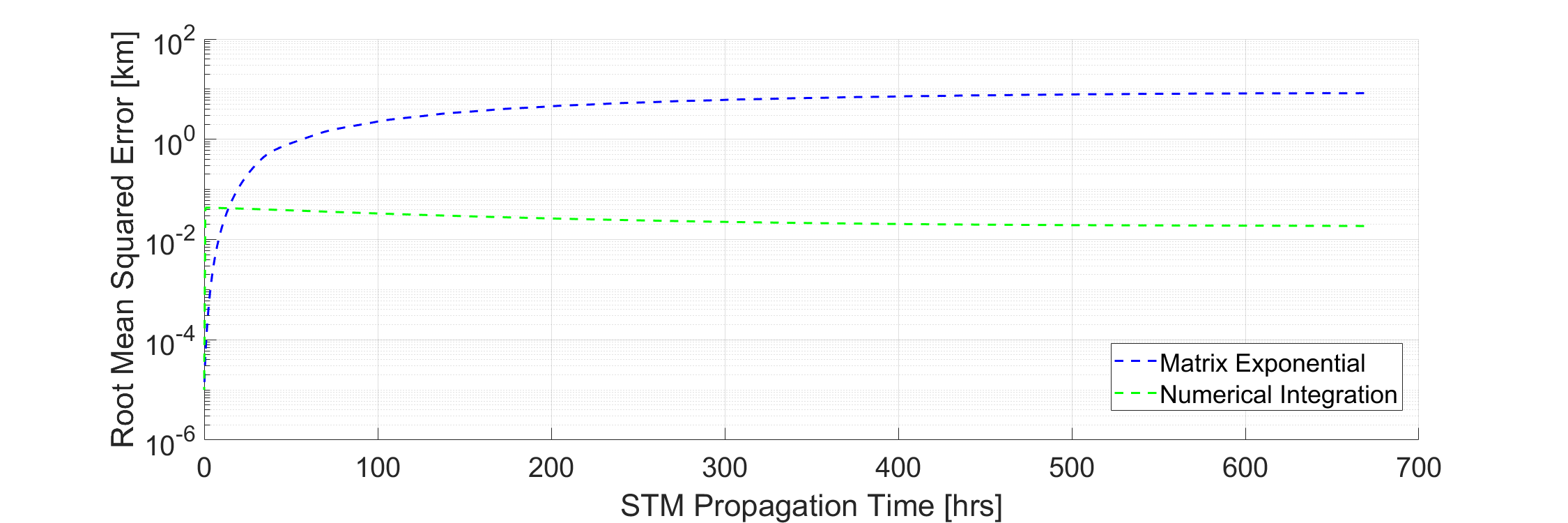}
    \caption{RMS propagation errors of each STM creation strategy from the initial chief and deputy states in Reconfiguration 1. }
    \label{fig: rms reconfig 1}
\end{figure}

\begin{figure}[h]
    \centering
    \includegraphics[width=1\linewidth]{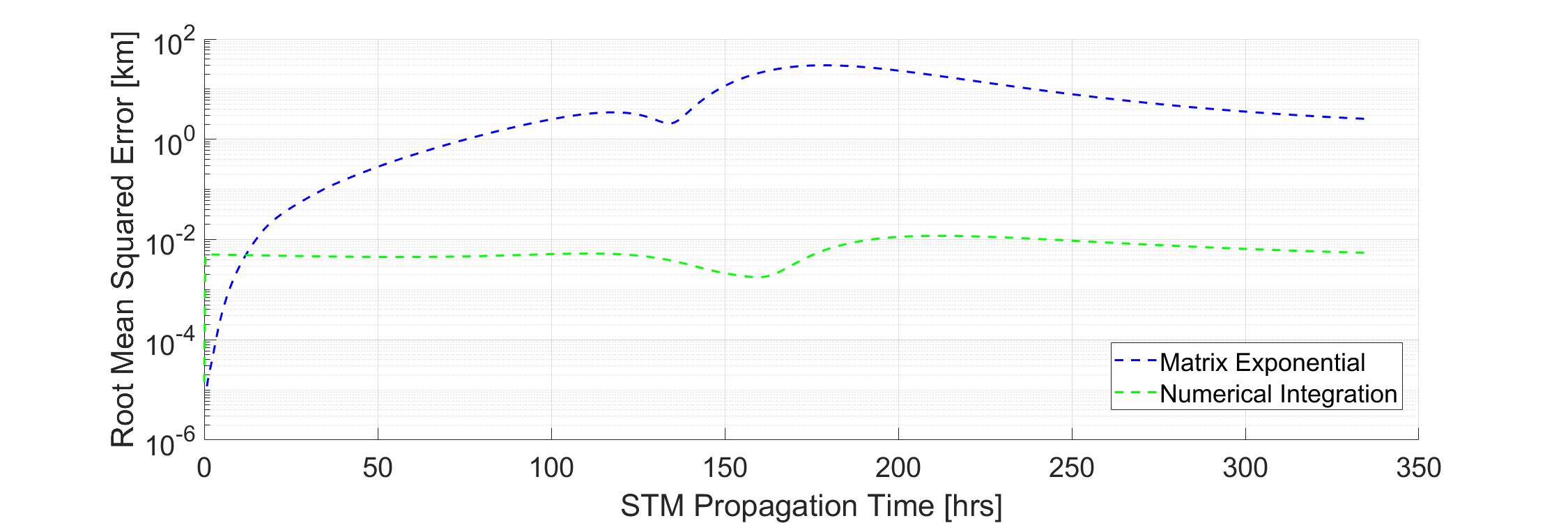}
    \caption{RMS propagation errors of each STM creation strategy from the initial chief and deputy states in Reconfiguration 2.}
    \label{fig: rms reconfig 2}
\end{figure}

As expected, the similar trajectories and maneuvers directions in Figure \ref{fig: reconfigs} correspond to similar optimal maneuver plans costs found by both STM strategies for each reconfiguration.
The STM and solver runtimes are both significantly longer for Reconfiguration 2 than Reconfiguration 1 despite having a far shorter control window.
This is due to the proximity of perilune, causing dynamic variability that requires more numerical effort to both calculate matrix exponentials and numerically integrate, as well as preventing quick solver convergence. 
The computational advantages of the matrix exponential approach are clear here, taking only 6.37\% and 51.05\% of the numerical integration approach for the two reconfigurations.
The two STM approaches also result in very similar solver runtimes, demonstrating that the choice of STM does not affect the computational performance of the optimization algorithm.
In contrast, the dynamic accuracy advantages of the numerical integration approach are clear in the final state terminal control error, 9.53\% and 1.66\% of the matrix exponential approach for Reconfiguration 1 and 2, respectively.
This difference in control accuracy can be further understood by examining the STM propagation errors resulting from the two STM approaches, shown in Figures \ref{fig: rms reconfig 1} and \ref{fig: rms reconfig 2} with y-axes on a log scale for ease of comparison.


The STM Root-Mean-Squared (RMS) propagation errors are found by propagating the initial deputy state with the chosen STM and comparing the resulting state against the numerically integrated, unperturbed CR3BP ground truth over the same time interval.
In both reconfigurations, the RMS propagation error of the numerical integration approach does not exceed $0.1$km while the RMS error of the matrix exponential approach can reach or exceed $10$km.
This difference in relative accuracy matches that for the final state RMS error for both reconfigurations in Table \ref{tab: reconfigs}.
Again, the proximity of perilune has a clear effect on performance. 
Both STM approaches in Reconfiguration 1 reach a fairly steady-state RMS error for propagation times greater than 100h, while passing over perilune in Reconfiguration 2 results in an oscillation in RMS error that only approaches similar steady-state values to Reconfiguration 1 for propagation times greater than 300h.
In general, the initial validation reconfigurations reflect the expected tradeoff in computational efficiency and dynamic accuracy present between the matrix exponential and numerical integration approaches.
Additionally, the accurate linear dynamic model resulting from the numerical integration strategy both near and far from perilune demonstrates the solver's capability to produce accurate maneuvers plans in both consistent and highly variable dynamic environments.

\subsection{Validation Against Two-Body Dynamics Models}

The Hill-Clohessy-Wiltshire (HCW) and Yamanaka-Ankersen (YA) dynamics models are commonly used to accurately characterize and control the relative motion of DSS in two-body Earth orbit. These models are used here as a baseline to evaluate the performance of the two CR3BP LTV relative motion modeling approaches.
The same problem setup from Reconfiguration 1 in the \textit{Initial Validation} section was used for this comparative testing, with the exception of the control window, set to 167.1h to evaluate the performance over longer durations, and the timestep, reduced to 1 minute in an attempt to increase the accuracy of the two-body dynamics models. As before, the solved maneuver plans are then simulated in a ground truth numerical propagation of the chief and deputy CR3BP unperturbed dynamics.

\begin{figure}[ht!]
    \centering
    \includegraphics[width=1\linewidth]{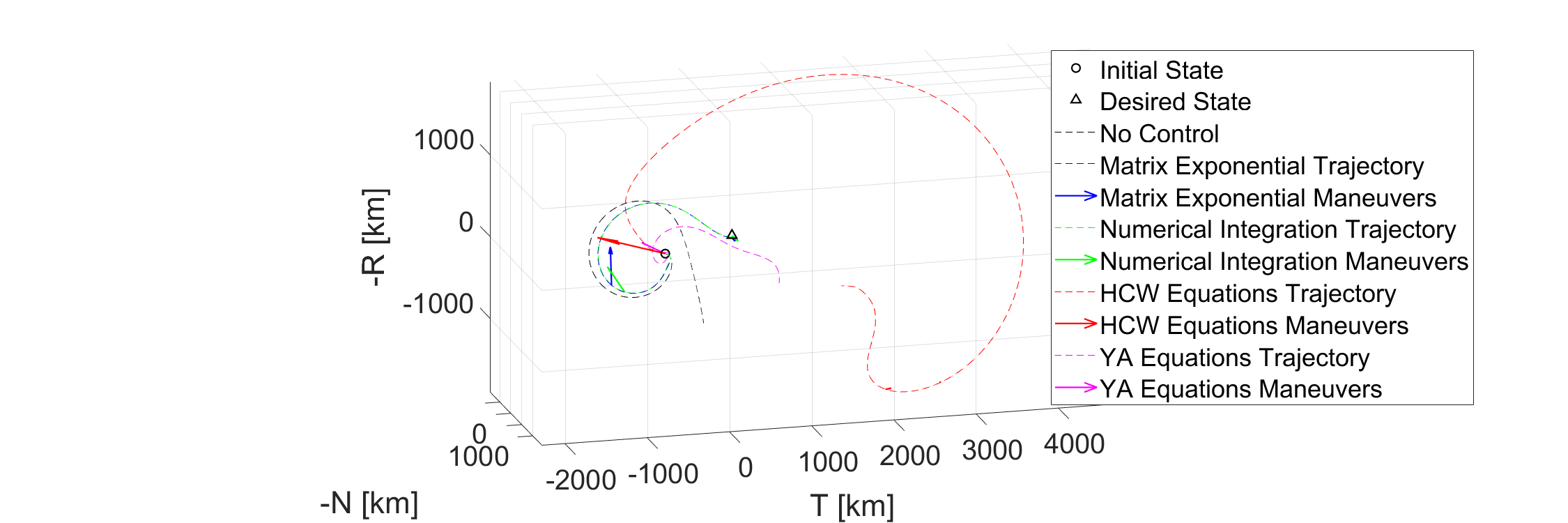}
    \caption{LVLH trajectories for Reconfiguration 1 with an extended control window.}
    \label{fig:unit-traj}
\end{figure}

\begin{figure}[h!]
    \centering
    \includegraphics[width=.82\linewidth]{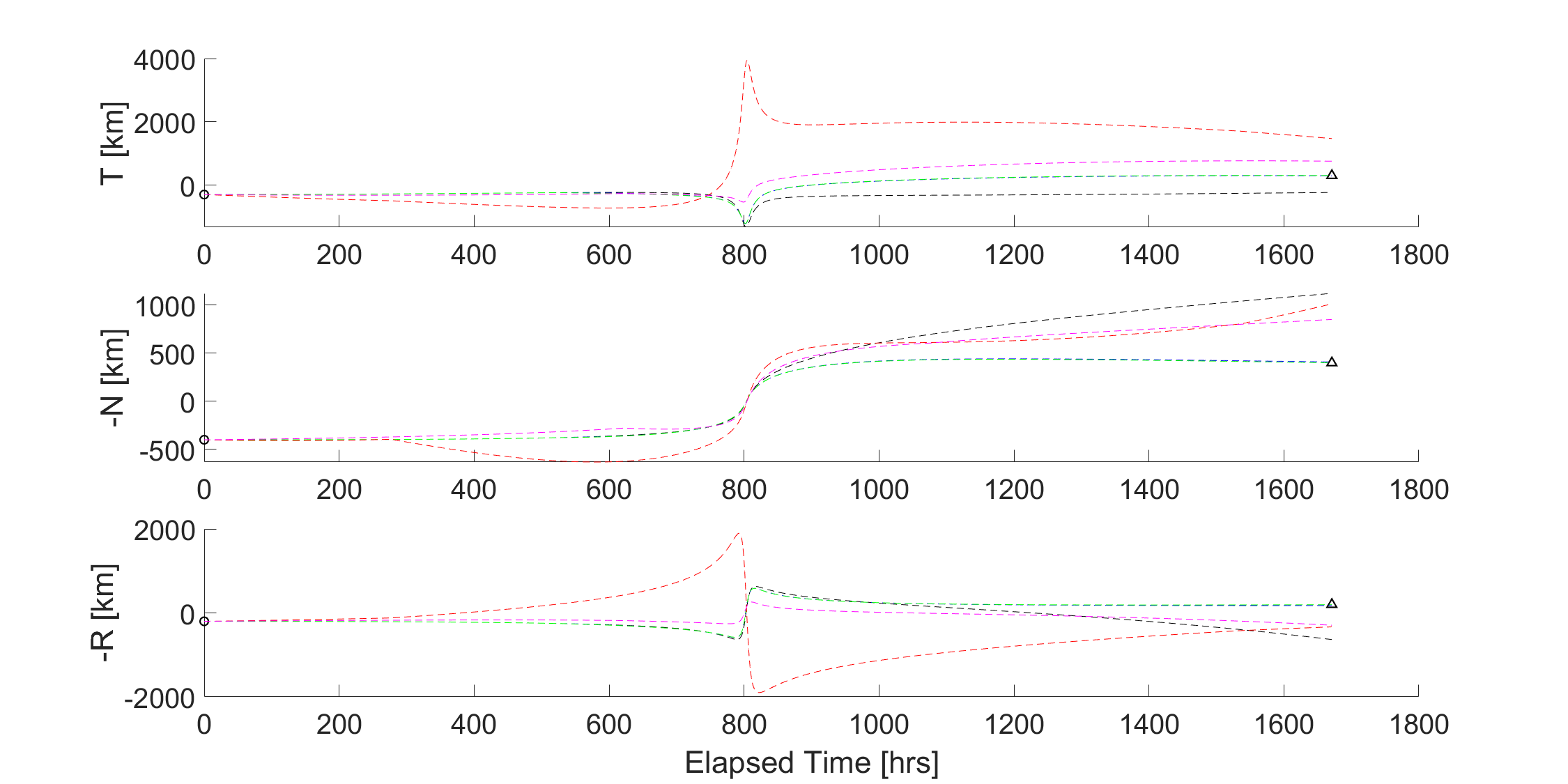}
    \caption{Time evolution of deputy LVLH position components for Reconfiguration 1 with an extended control window. Color coding is same as Figures \ref{fig:unit-traj} and \ref{fig:unit-err}.}
    \label{fig:unit-pos}
\end{figure}

\begin{figure}[h!]
    \centering
    \includegraphics[width=.82\linewidth]{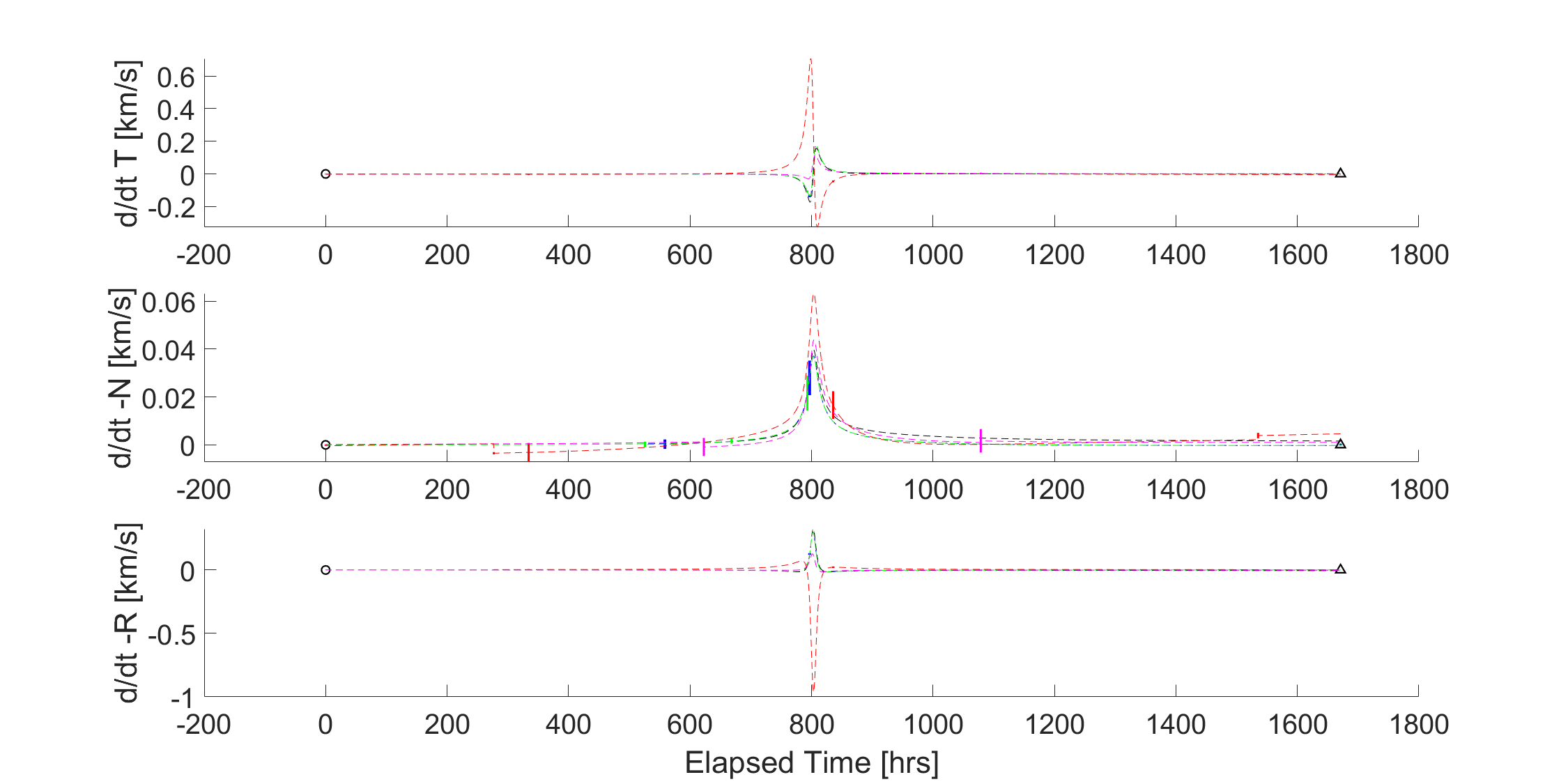}
    \caption{Time evolution of deputy LVLH velocity components for Reconfiguration 1 with an extended control window. Color coding is same as Figures \ref{fig:unit-traj} and \ref{fig:unit-err}.}
    \label{fig:unit-vel}
\end{figure}
\newpage

\begin{figure}[ht!]
    \centering
    \includegraphics[width=1\linewidth]{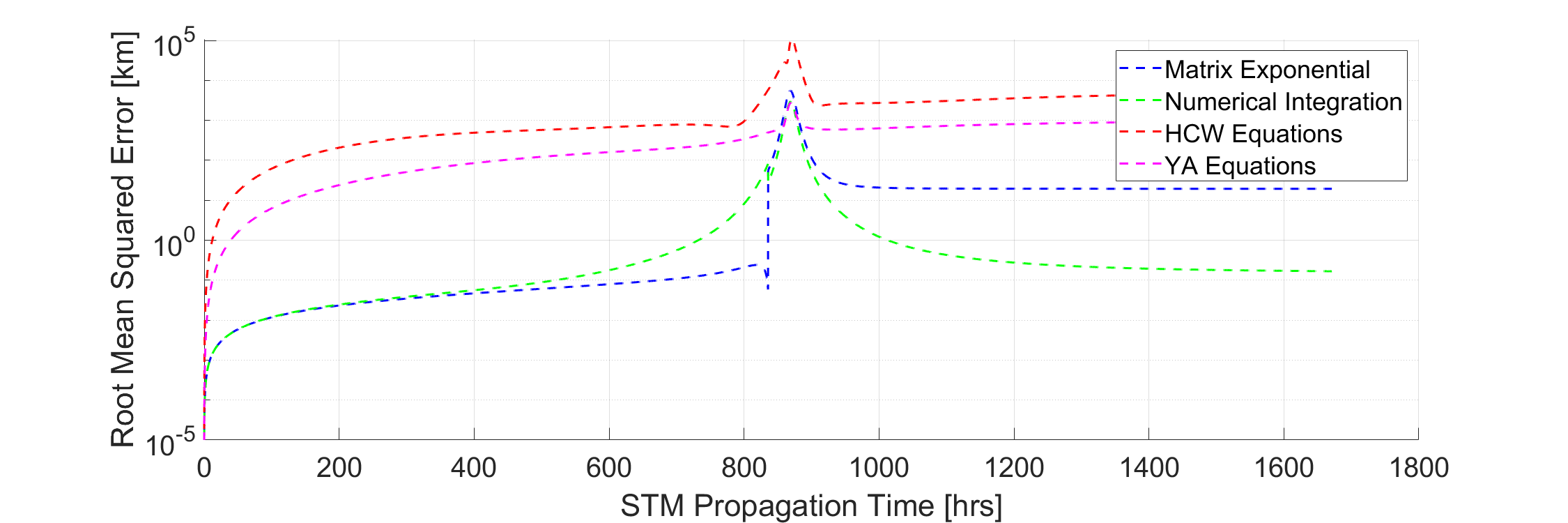}
    \caption{RMS propagation errors of each STM creation strategy from the initial chief and deputy states in Reconfiguration 1 with an extended control window.}
    \label{fig:unit-err}
\end{figure}

Figures \ref{fig:unit-traj}, \ref{fig:unit-pos}, and \ref{fig:unit-vel} display the trajectories of each dynamics model, where the STMs generated by each of the five models result in five distinctly different maneuver plans when input to the KD solver. Figure \ref{fig:unit-err} shows the final state propagation error of the different models, with the numerical integration approach yielding the lowest error, followed by the matrix exponential approach, while the two-body dynamics models all have drastically higher errors. A bifurcation between the CR3BP LTV and the two-body models can be seen almost immediately in the RMS error timeseries, and eventually this delta in error magnitude reaches over 100 km by the terminal state. Notably, all approaches show an increase in error when transiting perilune. Nevertheless, this example demonstrates the proposed CR3BP LTV models outperforming the baseline models, providing evidence of the benefits gained by employing dynamics models tailored to the CR3BP. 

\subsection{Monte Carlo Validation}

A Monte Carlo simulation is the next step taken in robustly validating the CR3BP LTV approaches over a wide range of possible scenarios. 
The goal of the Monte Carlo simulation is to 1) Demonstrate clear performance above the baseline two-body models (HCW and YA), 2) Quantify the tradeoffs in accuracy and computational efficiency between the numerical integration and matrix exponential approaches, and 3) Evaluate the solver performance over a wide range of cislunar reconfiguration scenarios. 
The setup of the Monte Carlo simulation, including which parameters are varied, and what distributions are used for sampling, is described in Table \ref{tab:montecarlo-setup}. 
Deputy initial state, deputy final state, chief initial absolute Halo orbit, and chief phase along the Halo orbit are all sampled from a representative distribution of different separations and scenarios.
The control window length is also chosen randomly to evaluate performance on short and long duration formation keeping or reconfiguration activities.
In total, 100 test cases are sampled and used in the simulation. All five dynamics models are evaluated on the same 100 test cases for consistency. The time step used for the LTI approximation within the matrix exponential approach is set at 60s for all test cases, and the KD solver setup matches that used in \textit{Initial Validation}.

\begin{table}[h!]
    \centering
    \begin{tabular}{c|c}
         Variable  &  Sampling Approach\\
         \hline\hline
         Chief initial absolute state & Uniform sampling from 6000 known halo orbit states: \\
         & 1000 states each from the 9:2, 4:1, 7:2, 3:1, 5:2, 2:1 halo orbits \\
         Deputy initial relative position & Log-uniform distribution $\sim$ [-5000, 5000] km\\
         Deputy initial relative velocity & Normal distribution, Mean = 0 km/s, Std. dev. = 0.001 km/s\\
         Deputy final relative position & Log-uniform distribution $\sim$ [-5000, 5000] km\\
         Deputy final relative velocity & Normal distribution, Mean = 0 km/s, Std. dev. = 0.001 km/s\\
         Control window   & Log-uniform distribution $\sim$ [0.1$\pi$, 4$\pi$] TU\\
         \hline
    \end{tabular}
    \caption{Variables of interest and associated sampling approaches for Monte Carlo simulation.}
    \label{tab:montecarlo-setup}
\end{table}

\begin{figure}[h!]
    \centering
    \includegraphics[width=0.6\linewidth]{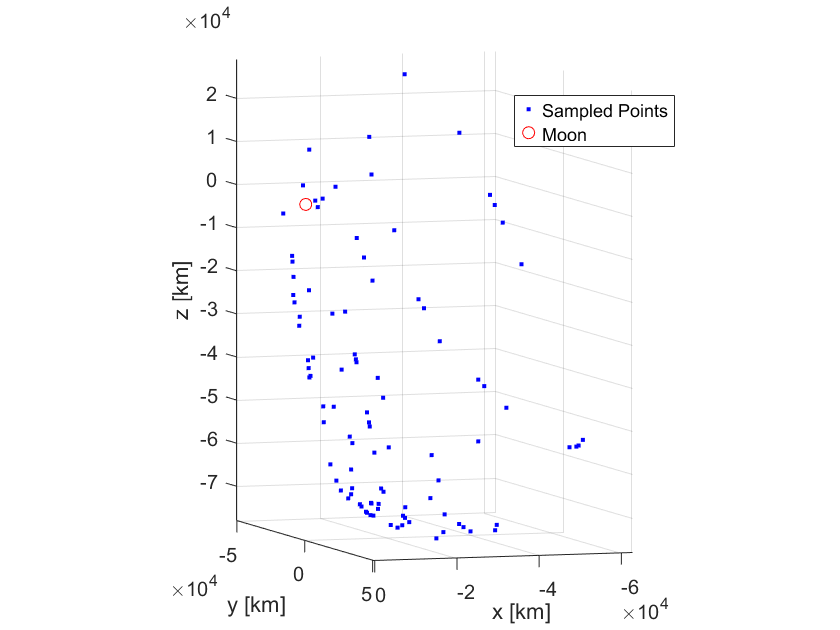}
    \caption{Chief spacecraft initial absolute states (in moon-centered synodic frame) used in Monte Carlo simulation, sampled randomly from 9:2, 4:1, 7:2, 3:1, 5:2, and 2:1 Halo orbits.}
    \label{fig:montecarlo-x0}
\end{figure}

Figure \ref{fig:montecarlo-x0} provides a visualization of the sampling approach taken for the chief spacecraft initial absolute state, drawing from multiple halo orbits and different phases along those orbits, including at perilune.

\begin{table}[h!]
\setlength{\tabcolsep}{5pt}
    \centering
    \begin{tabular}{l|c|c|c|c}
        \multicolumn{1}{c|}{Metric} & Matrix Exp. & Num. Int. & HCW & YA \\
        \hline \hline
        Final Position Error (km) - Median & 8.8225 & 9.8411 & 6208.8 & 3855.3\\
        Final Position Error (km) - Min & 7.2799e-2 & 7.2862e-2 & 26.469 & 27.665\\
        Final Position Error (\%) - Median & 5.1399 & 3.5769 & 2620.9 & 1163.9\\
        Final Position Error (\%) - Min & 7.4492e-3 & 2.9771e-3 & 17.707 & 9.7436\\
        Optimal Cost (m/s) - Median & 14.547 & 14.542 & 10.265 & 5.3085\\
        Optimal Cost (m/s) - Mean & 35.166 & 35.661 & 18.812 & 20.290\\
        Optimal Cost (m/s) - Max & 344.15 & 348.67 & 97.508 & 351.65\\
        Optimal Cost (m/s) - Min & 1.8811 & 1.4988 & 0.98025 & 0.57857\\
        STM Runtime (s) - Median & 1.9914 & 12.597 & 1.5156 & 1.8405\\
        STM Runtime (s) - Mean & 4.9494 & 22.369 & 3.4174 & 3.9759\\
        STM Runtime (s) - Max & 23.727 & 84.492 & 16.541 & 16.864\\
        STM Runtime (s) - Min & 0.19371 & 1.6645 & 0.13326 & 0.17465\\
        Solver Runtime (s) - Median & 7.06665 & 7.4127 & 6.7330 & 6.9115\\
        Solver Runtime (s) - Mean & 20.225 & 16.626 & 15.865 & 22.378\\
        Solver Runtime (s) - Max & 285.08 & 87.209 & 94.789 & 185.90\\
        Solver Runtime (s) - Min & 1.0622 & 1.0237 & 0.68229 & 0.98903\\
        \hline
    \end{tabular}
    \caption{Performance metrics for Monte Carlo simulation, comparing multiple dynamics models.}
    \label{tab:montecarlo}
\end{table}

Table \ref{tab:montecarlo} reports the results of the Monte Carlo simulation. 
Clearly, at scale, the same trends persist from the individual test case comparing the five dynamics models. The matrix exponential and numerical integration approaches had a median terminal position error of just 5.1\% and 3.6\%, respectively. 
In contrast, the HCW and YA models failed to model the dynamics accurately enough for the solver to produce any accurate solutions, with median terminal position error of over 1000\%.
When comparing the matrix exponential and numerical integration approaches, the data suggests that the tradeoffs in accuracy and efficiency favor the matrix exponential approach. 
The average time needed for STM generation is an order of magnitude lower for the matrix exponential approach compared to numerical integration, while the average solver runtimes are comparable. However, the matrix exponential approach does not suffer a significant drop in error compared to numerical integration. 
In short, the matrix exponential approach appears to be much faster than numerical integration, while only sacrificing a relatively small level of accuracy on average.

\subsection{Model Predictive Control}

The final component of analysis and validation for this work is to implement a Model Predictive Control (MPC) scheme. 
MPC rejects dynamic propagation errors, which stem from the inherent differences between the matrix exponential and the numerically integrated LTV model, and can also mitigate maneuver execution errors and navigation uncertainty. 
The naive MPC architecture begins with an initial control solution at the beginning of the control window and periodically produces a new maneuver plan for the remaining control window from the current state estimate.
This strategy aims to leverage the more computationally efficient matrix exponential approach to enable accurate, robust closed loop control.
For the MPC test case, the initial and desired final states of the deputy are given in Table \ref{tab:mpc-initial and final deputy states}, and the control window was 167.1h, to simulate a long-term reconfiguration. 
The MPC strategy divides this control window into 10 segments and then re-solves at the end of each segment, as the deputy spacecraft draws closer to the terminal state. 
The chief spacecraft initial absolute state is at apolune of the 9:2 NRHO.
Typical maneuver execution and navigation errors for small satellites in the CR3BP regime are considered as part of the MPC validation simulation environment.
Navigation and maneuver execution performance are simulated by a set of zero-mean Gaussian distributions given in Table \ref{tab:mpc-cov}.

\begin{table}[h!]
    \centering
    \begin{tabular}{c|c|c|c|c|c|c}
    Deputy State  & $x$ (km) & $y$ (km) & $z$ (km) & $\dot{x}$ (km/s) & $\dot{y}$ (km/s) & $\dot{z}$ (km/s) \\
    \hline \hline
    Initial & -3000 & -4000 & -2000 & 0 & 0 & 0 \\
    Final &  3000 & 4000 & 2000 & 0 & 0 & 0 \\
    \hline
    \end{tabular}
    \caption{Initial and final deputy relative states in LVLH frame for MPC reconfiguration.}
    \label{tab:mpc-initial and final deputy states}
\end{table}

\begin{table}[h!]
    \centering
    \begin{tabular}{c|c}
         Parameter &  Covariance \\
         \hline\hline
         Chief absolute position components & $1^2$ km$^2$ \\
         Chief absolute velocity components & $0.01^2$ km$^2$/s$^2$\\
         Deputy relative position components & $0.01^2$ km$^2$\\
         Deputy relative velocity components & $0.001^2$ km$^2$/s$^2$\\
         Maneuver time & $60^2$ s$^2$\\
         Maneuver magnitude & $0.01^2$ km$^2$\\
         Maneuver direction & $1^2$ deg$^2$\\
         \hline
    \end{tabular}
    \caption{Navigation and maneuver execution uncertainties used in MPC reconfiguration.}
    \label{tab:mpc-cov}
\end{table}

\begin{figure}[ht!]
    \centering
    \includegraphics[width=1\linewidth]{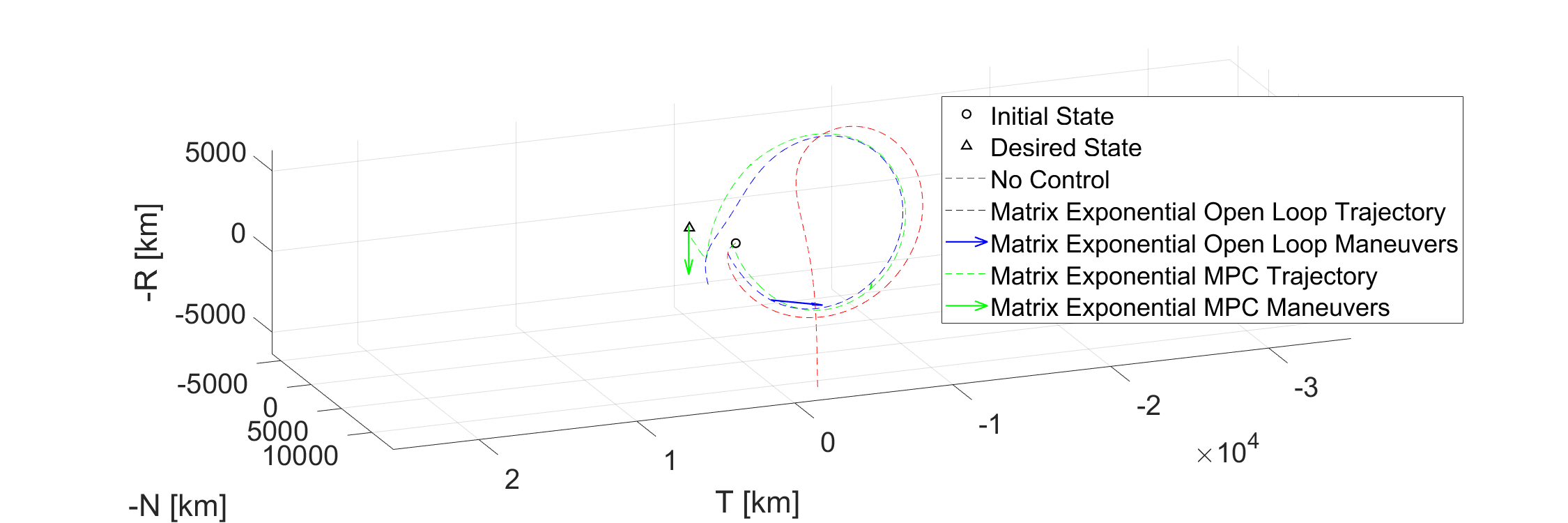}
    \caption{Trajectory of deputy in LVLH frame, comparing MPC to open loop control.}
    \label{fig:mpc-traj}
\end{figure}

\begin{figure}[h!]
    \centering
    \includegraphics[width=.82\linewidth]{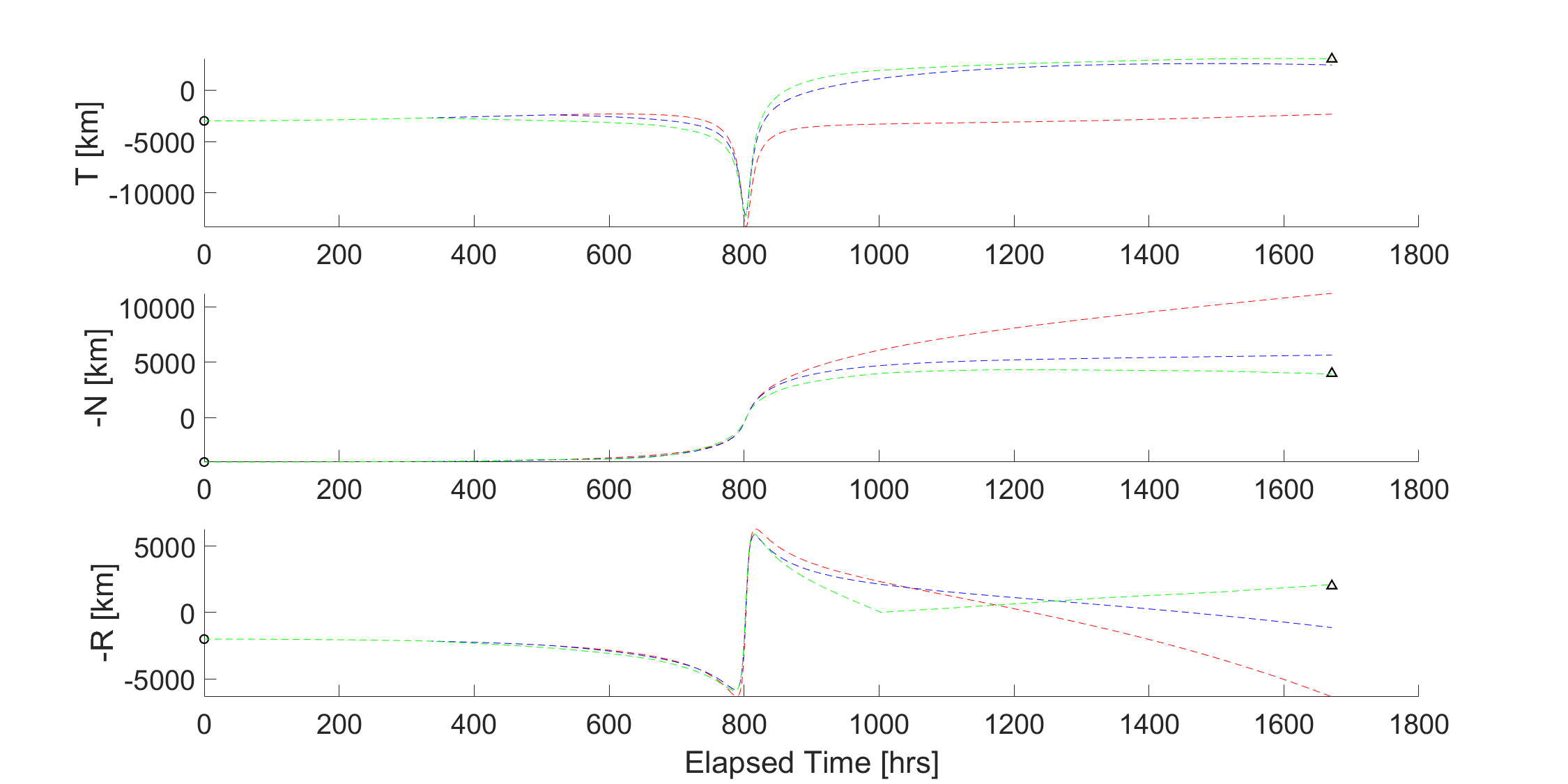}
    \caption{Time evolution of deputy position components in LVLH frame, comparing MPC to open loop control.}
    \label{fig:mpc-pos}
\end{figure}

\begin{figure}[h!]
    \centering
    \includegraphics[width=.82\linewidth]{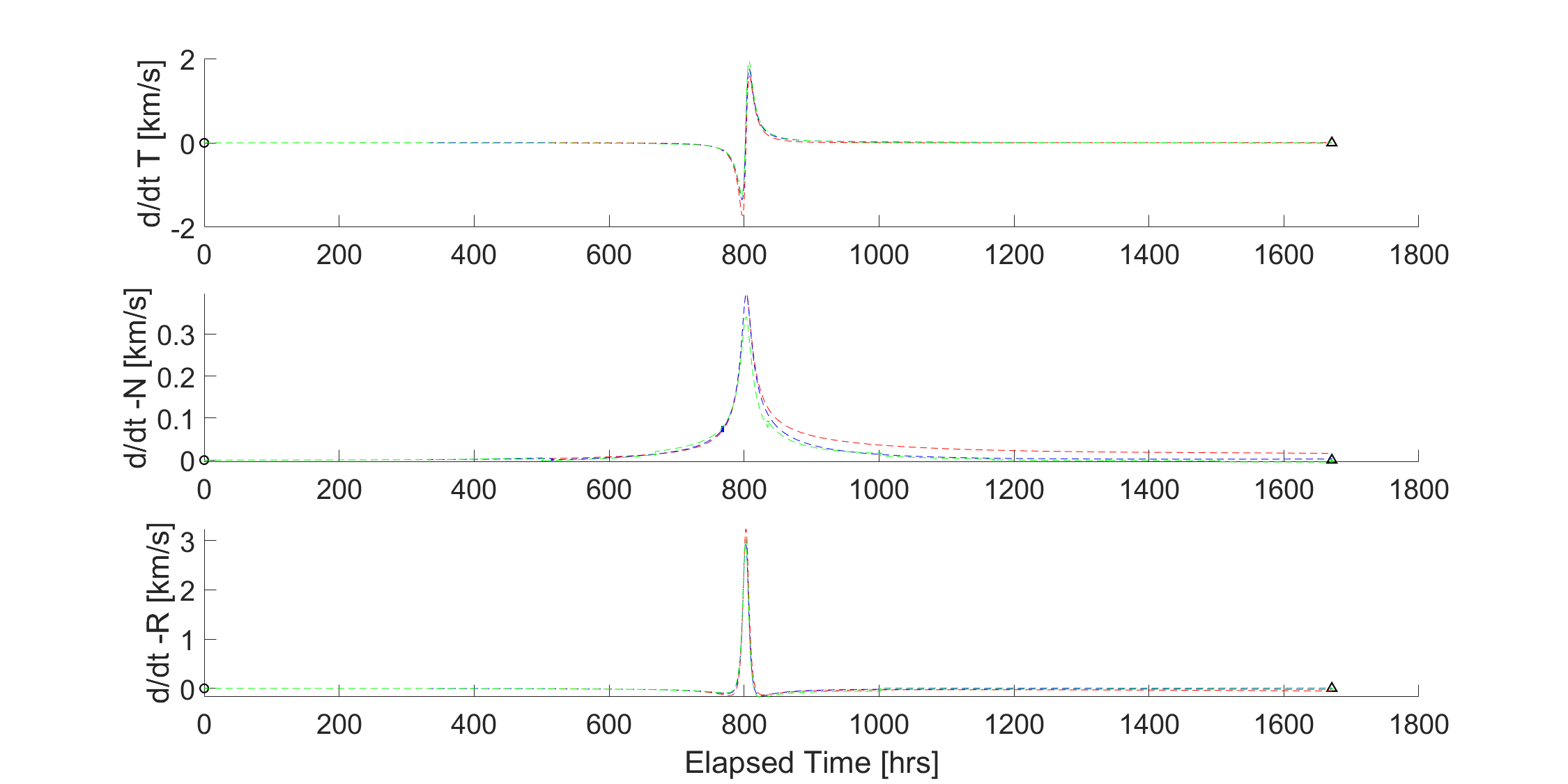}
    \caption{Time evolution of deputy velocity components in LVLH frame, comparing MPC to open loop control.}
    \label{fig:mpc-vel}
\end{figure}

The results of the MPC simulation for a single test case can be visualized in Figures \ref{fig:mpc-traj}, \ref{fig:mpc-pos}, and \ref{fig:mpc-vel}. 
Open loop control and uncontrolled trajectories are included for comparison. 
Open loop control in this context is equivalent to the naive maneuver planning used in previous sections of this paper. 
The tradeoffs between MPC and open loop control are summarized in \ref{tab:mpc-metrics}, and are consistent with expectations. The terminal position error from MPC is less than a tenth of that attained using open loop control, due to the MPC regularly re-planning as the time horizon shrinks. However, this increase in accuracy comes at the expense of $\Delta v$ cost and runtime, where the metrics for MPC exceed those for open loop control by a factor of about four.

\begin{table}[h!]
    \centering
    \begin{tabular}{l|c|c}
        \multicolumn{1}{c|}{Metric} & MPC & Open Loop \\
        \hline \hline
        Terminal Position Error (km) & 162.25 & 2087.0 \\
        Terminal Position Error (\%) & 3.0130 & 38.755 \\
        Optimal Cost (m/s) & 92.885 & 23.131 \\
        \hline
    \end{tabular}
    \caption{Performance metrics for MPC validation reconfigurations.}
    \label{tab:mpc-metrics}
\end{table}

\section{Conclusions}

This work applies a generic reachable set theory optimal control approach for Linear Time-Variant (LTV) dynamic systems to cislunar relative spacecraft motion. 
First, a linear model of Circular Restricted Three-Body Problem (CR3BP) relative spacecraft motion is leveraged to create an LTV dynamics model. 
The instantaneous changes in relative Cartesian state are numerically integrated to produce a State Transition Matrix (STM) that translates impulsive control actions from time of execution to the end of the control window. 
Second, the LTV model is incorporated into a reachable set theory-based solver that uses a set of optimality conditions to efficiently produce impulsive control profiles with provable optimality.
Accommodations are made within the optimization algorithm to make it robust to the extreme variations in relative and absolute dynamics present near perilune. 
Third, this methodology demonstrates accurate control performance over short and long control windows, accomplishing both small and large reconfigurations over different CR3BP orbits. 
These capabilities are enhanced by a Model Predictive Control (MPC) architecture to reject all sources of control, navigation, and dynamic error.

The MPC results are promising, as they demonstrate the capability to handle uncertainties and perturbations through closed loop control. 
While the work presented in this paper is subject to CR3BP dynamics, it could be extended to higher fidelity models in the future. 
Adding perturbations such as SRP and lunar spherical harmonics to the ground truth simulation is a logical next step, with the ultimate goal of using a full-ephemeris model. 
Other higher fidelity models, such as the Elliptical Restricted Three-Body Problem (ER3BP) could be employed as well. 
It is important to note that since an LTV system is required to be compatible with the KD solver, there will always be a delta between the STM dynamics and the ground truth dynamics due to the necessary linearizations. 
Further investigation is required to determine whether the KD solver is capable of overcoming this delta in higher fidelity models, both in open loop and closed loop control approaches.
Overall, this work demonstrates an optimal control methodology for cislunar applications with the computational efficiency to be used on-board spacecraft, enabling spacecraft autonomy beyond Earth orbit.


\section{Acknowledgments}
This work is supported by Blue Origin (SPO \#299266) as an Associate Member and Co-Founder of the Stanford’s Center of AEroSpace Autonomy Research (CAESAR), and by the National Science Foundation Graduate Research Fellowship Program under Grant No. DGE-2146755. This article solely reflects the opinions and conclusions of its authors and not any of its sponsors.

\appendix

\bibliographystyle{AAS_publication}   
\bibliography{references}   

\end{document}